\documentclass[11pt]{article}

\usepackage{amsmath}
\usepackage{amsfonts}
\usepackage[portuges]{babel}
\usepackage{amssymb}
\usepackage{latexsym}
\usepackage{enumerate}
\usepackage[dvips]{graphicx}

\def\rain{\to \infty}

\def\N{{\rm I\kern-.20em N}}
\def\R{{\rm I\kern-.20em R}}
\def\indi{{1\kern-.20em\rm I}}

\def\refhg{\hangindent=15pt\hangafter=1}
\def\refmark{\par\vskip 2mm\noindent\refhg}
\def\bkR{{\rm I\kern-.17em R}}
\def\bkN{{\rm I\kern-.20em N}}

\newtheorem{lema}{Lemma}[section]

\newtheorem{coro}{Corollary}[section]
\newtheorem{prop}{Proposition}[section]

\newtheorem{defi}{Definition}[section]
\newtheorem{ex}{Example}[section]

\newcommand{\pg}{\hspace{0.6cm}}

\newcommand{\bdem} {\begin{proof}}
\newcommand {\edem}{\hfill $\square$ \end {proof}}

\begin{document}
  \title{On the multivariate upcrossings index}
\author{C. Viseu \footnote{Instituto Superior de Contabilidade e
Administra\c c\~ ao de Coimbra, Coimbra,
Portugal.\newline
E-mail:cviseu@iscac.pt}\\  Instituto Superior de Contabilidade\\
e Administra\c c\~ ao de Coimbra\\ Portugal \and L. Pereira, A.P.
Martins, H. Ferreira \footnote{Departamento de
Matem\'atica, Universidade da Beira Interior, 6200 Covilh\~a,
Portugal.
E-mail: lpereira@ubi.pt; amartins@ubi.pt;
helena.ferreira@ubi.pt}\\Departamento de Matem\'atica \\
Universidade da Beira Interior\\ Portugal}
\date{}
\maketitle

\noindent {\bf Abstract:} The notion of multivariate upcrossings
index of a stationary sequence
${\bf{X}}=\{(X_{n,1},\ldots,X_{n,d})\}_{n\geq 1}$ is introduced and
its main properties are derived, namely the relations with the
multivariate extremal index and the clustering of upcrossings.

Under asymptotic independence conditions on
the marginal sequences of ${\bf{X}}$ the multivariate upcrossings
index is obtained from the marginal upcrossings indices.

For a class of stationary multidimensional sequences assumed to satisfy a mild oscillation res\-tric\-tion, the
multivariate upcrossings index is computed from the joint
distribution of a finite number of variables.

The upcrossings index is calculated for two examples of bivariate
sequences.\vspace{0.5cm}


\section{Introduction}

\pg Let ${\bf{X}}=\{{\bf{X}}_n=(X_{n,1},\ldots,X_{n,d})\}_{n\geq 1}$
be a $d-$dimensional stationary sequence. We will denote the
variable  $\max\{X_{i,j},\ i\in I\}$ by $M_{n,j}(I),$ $I\subset \N,$
$j=1,\ldots,d,$ and consider
${\bf{M}}_n(I)=(M_{n,1}(I),\ldots,M_{n,d}(I)).$ When
$I=\{1,\ldots,s\},$  we shall simply write ${\bf{M}}_s$ and
$M_{s,j}.$

For ${\bf{a}}=(a_1,\ldots,a_d),$  ${\bf{b}}=(b_1,\ldots,b_d)$ and
${\bf{x}}=(x_1,\ldots,x_d)$ in $\R^d$ and $c\in \R,$ let
${\bf{ax+b}}=$ $(a_1x_1+cb_1,\ldots,a_dx_d+cb_d),$
${\bf{a}}^c=(a_1^c,\ldots,a_d^c)$ and ${\bf{a\leq b}}$ if and only
if $a_j\leq b_j$ for all $j=1,\ldots,d.$\vspace{0.5cm}

We shall assume that the $d$-dimensional stationary sequence
${\bf{X}}$ satisfies the mixing condition $\Delta({\bf{u}})$ of
Nan\-da\-go\-pa\-lan (1990), which we next recall.

\begin{defi}
The sequence ${\bf{X}}$ is said to satisfy the condition
$\Delta({\bf{u}})$, for some sequence
${\bf{u}}=\{{\bf{u}}_n=(u_{n,1},\ldots,u_{n,d})\}_{n\geq 1}$ of
elements of $\R^d,$ if $\alpha_{n,l_{n}}\to 0$ as $n\to \infty $ for
some sequence $l_{n}=o(n),$ where
$$\alpha _{n,l}=sup\left\{ \left| P(A\cap B)-P(A)P(B)\right|:
\ A\in \mathcal{B}_1^k({\bf{u}}),\ B\in
\mathcal{B}_{k+l}^n({\bf{u}}),\ 1\leq k\leq n-l\right\}
$$ and $B_i^j({\bf{u}})$ denotes the $\sigma$-field generated by the events
$\left\{ X_{s,k}\leq u_{n,k}\right\},$ $k=1,\ldots,d,$ $i\leq s\leq
j.$
\end{defi}

Condition $\Delta({\bf{u}})$ implies condition $D({\bf{u}})$ of
Hsing (1989) and, for ${\bf{u}}_j=\{u_{n,j}\}_{n\geq 1},$ condition
$\Delta({\bf{u}}_j)$ for ${\bf{X}}_j=\{X_{n,j}\}_{n\geq 1},$
$j=1,\ldots,d.$\vspace{0.5cm}

In order to study the asymptotic behaviour of the upcrossings
$\{X_{i,j}\leq u_{n,j}<X_{i+1,j}\},$ $j=1,\ldots,d,$ among the first
$n$ variables of ${\bf{X}},$ we consider the vector of point
processes
${\bf{S}}_n^{({\bf{u}})}=(S_{n,1}^{(u_{n,1})},\ldots,S_{n,d}^{(u_{n,d})}),$
$n\geq 1,$ with
$$S_{n,j}^{(u_{n,j})}(B)=\sum_{i=1}^n\indi_{\{X_{i,j}\leq
u_{n,j}<X_{i+1,j}\}}\delta_{\frac{i}{n}}(B),\quad B\subset [0,1],\
j=1,\ldots,d,$$ where $\indi_{A}$ denotes the indicator of event $A$
and $\delta _{a}$ the unit mass at $a.$ Thus, $S_{n,j}^{(u_{n,j})}$
has unit mass at $i/n$ if $X_{i,j}\leq u_{n,j}<X_{i+1,j},$
$j=1,\ldots,d.$ We shall omit, in the definition of the point
processes, the reference to the sequence of levels ${\bf{u}},$ which
defines the upcrossings, whenever it is clear.\vspace{0.3cm}

For a broad class of weakly dependent sequences, Ferreira (2006)
proved that the limit in distribution of $S_{n,j},$ $j=1,\ldots,d,$
if it exists, is a compound Poisson process. When the limiting mean
number of upcrossings of the level $u_{n,j}$ is $\nu_{j}>0,$ the
Poisson rate of the limiting point process is $\eta_j\nu_j,$ being
$\eta_j$ the uprcossings index of the stationary sequence
${\bf{X}}_j.$ The upcrossings index is related with the extremal
index $\theta_j$ through $\eta_j=\tau_j\theta_j/\nu_j,$ where
$\tau_j$ is the limiting mean number of exceedances of the level
$u_{n,j}.$ Ferreira (2007) also proved that the reciprocal of the
upcrossings index can be interpreted as the limiting mean size of
clusters of upcrossings.

In this paper we extend the univariate theory described above to the
multivariate case.

In Section 2 we study the vector the upcrossing point processes and
we obtain some weak convergence results.

The multivariate upcrossings index is then defined in Section 3 and
its properties and relations with the multivariate extremal index
and its marginal upcrossings indices are established.

Under asymptotic independence conditions on the marginal sequences
of ${\bf{X}},$ in Section 4, we obtain the multivariate upcrossings
index from the marginal indices.

In Section 5, by imposing local restrictions on the oscillations of
${\bf{X}},$ we compute the multivariate upcrossings index from a
finite number of the variables of ${\bf{X}}.$

Finally in Section 6, two examples of bivariate sequences are
presented for which the upcrossings index is calculated.

\section{Upcrossing point processes}

\pg Under condition $\Delta({\bf{u}}),$ Nandagopalan (1990)
established a necessary and sufficient condition for the convergence
in distribution of the sequence of $d$-dimensional processes of
exceedences $\{X_j>u_{n,j}\},$ $j=1,\ldots,d,$ generated by a
sequence ${\bf{X}},$ not necessarily stationary. The arguments used
can easily be adapted to the $d$-dimensional sequence
$\{{\bf{S}}_n\}_{n\geq 1}$ of point processes of upcrossings. We
obtain in this way an extension of the theory presented in Ferreira
(2006), with the introduction of the concept of a multivariate
upcrossings index and its relations with the multivariate extremal
index and corespondent marginal indices.\vspace{0.5cm}

Let $\tilde{{\bf{u}}}^{({\boldsymbol{\nu}})},$
${\boldsymbol{\nu}}=(\nu_1,\ldots,\nu_d)\in \R^d_+,$ denote the
levels ${\bf{u}}$ such that
\begin{equation}
nP(X_{i,j}\leq u_{n,j}<X_{i+1,j})\xrightarrow [n\rain]{} \nu_j,\quad
j=1,\ldots,d,\label{eq2_1}
\end{equation}
and $\tilde{{\bf{u}}}_j^{(\nu_j)}$ denote the corresponding marginal
sequences ${\bf{u}}_j=\{u_{n,j}\}.$

Analogously, levels ${\bf{u}}$ such that
\begin{equation}
nP(X_{i,j}>u_{n,j})\xrightarrow [n\rain]{} \tau_j,\quad
j=1,\ldots,d,\label{eq2_2}
\end{equation}
will be denoted by ${\bf{u}}^{({\boldsymbol{\tau}})},$
${\boldsymbol{\tau}}=(\tau_1,\ldots,\tau_d)\in \R^d_+,$ and the
corresponding marginal sequences by
${\bf{u}}_j^{(\tau_j)}.$\vspace{0.5cm}

Note that if ${\bf{u}}=\tilde{{\bf{u}}}^{({\boldsymbol{\nu}})}$ and
${\bf{k}}=\{k_n\}_{n\geq 1}$ is such that $k_nl_n/n\xrightarrow
[n\rain]{} 0$ then, for any interval $I\subset [0,1]$ with length
$\leq l_n/n,$ it holds
\begin{equation}
k_nP({\bf{S}}_n(I)\neq {\bf{0}})\xrightarrow [n\rain]{}
0.\label{eq2_3}
\end{equation}

Let ${\bf{k}}=\{k_n\}_{n\geq 1}$ be a sequence of integers
satisfying
\begin{equation}
k_n\xrightarrow [n\rain]{}\infty ,\qquad \frac{k_nl_n}{n}
\xrightarrow [n\rain]{}0, \qquad k_n\alpha_{n,l_n}\xrightarrow
[n\rain]{} 0,\label{eq2_4}
\end{equation}
and $r_n=[n/k_n],\ n\geq 1.$ For such sequences the upcrossings
$\{X_{i,j}\leq u_{n,j}<X_{i+1,j}\}$ with  $i\in
\{(s-1)r_n+1,\ldots, sr_n\},$ for some $s=1,\ldots,k_n,$ are said to
belong to the same cluster of upcrossings of the $j-$th marginal
sequence.\vspace{0.5cm}

We establish, in the next lemma, the asymptotic independence of
upcrossings over disjoint intervals, analogously to Lemma 3.3.2 in
Nandagopalan (1990).

\begin{lema}
Assume condition $\Delta({\bf{u}})$ holds for ${\bf{X}}$ and
${\bf{k}}$ satisfies (\ref{eq2_3}) and (\ref{eq2_4}). Then, for any
${\bf{y}}_1,\ldots,{\bf{y}}_{k_n}\in \N_0^d$ and
$I_{n,1},\ldots,I_{n,k_n}$ disjoint subintervals of $[0,1]$, we have
$$P\left(\bigcap_{j=1}^{k_n}\{{\bf{S}}_n(I_{n,j})={\bf{y}}_j\}\right)-
\prod_{j=1}^{k_n}P({\bf{S}}_n(I_{n,j})={\bf{y}}_j)\xrightarrow
[n\rain]{} 0.$$
\end{lema}\vspace{0.5cm}

This lemma extends Lemma 2.2 of Hsing {\it et al.} (1988) to the
$d-$dimensional case and it's essential to obtain the following
necessary and sufficient condition for the convergence of
$\{{\bf{S}}_n\}_{n\geq 1}$ to a $d-$dimensional compound Poisson
process ${\bf{S}}$ with Laplace transform given by
$$L_{{\bf{S}}}(f_1,\ldots,f_d)=\exp\left(-\alpha\int_{[0,1]}\int_{\N_0^d-\{{\bf{0}}\}}
\left(1-\exp\left(-\sum_{j=1}^d y_jf_j(x)\right)\right)
d\Pi({\bf{y}})dx\right),$$ for each nonnegative, measurable
functions $f_1,\ldots,f_d$ on $[0,1],$ where $\alpha$ is a positive
constant and $\Pi$ a probability distribution on
$\N_0^d-\{{\bf{0}}\}.$ We shall represent such a point process by
${\bf{S}}[\alpha,\Pi]$ to make the dependence on its intensity
$\alpha$ and multiplicity distribution $\Pi$ explicit.

Let us consider the following partition of $[0,1]:$
$C_{n,1}=\left[0,\frac{r_n}{n}\right],$
$C_{n,j}=\left((j-1)\frac{r_n}{n},j\frac{r_n}{n}\right],$
$j=1,\ldots,k_n,$
$C_{n,k_n+1}=\left(k_n\frac{r_n}{n},1\right].$\vspace{0.2cm}

\begin{prop}
Assume condition $\Delta({\bf{u}})$ holds for ${\bf{X}}$ and
${\bf{k}}$ satisfies (\ref{eq2_3}) and (\ref{eq2_4}). Then
${\bf{S}}_n\xrightarrow [n\rain]{} {\bf{S}}[\alpha,\Pi]$ if and only
if $$P({\bf{S}}_n([0,1])={\bf{0}})\xrightarrow [n\rain]{}
e^{-\alpha}$$ and $$\Pi_n({\bf{y}})\equiv
P({\bf{S}}_n(C_{n,1})={\bf{y}}\ |\ {\bf{S}}_n(C_{n,1})\neq
{\bf{0}})\xrightarrow [n\rain]{} \Pi({\bf{y}}),\quad \forall
{\bf{y}} \in \N_0^d-\{{\bf{0}}\}.$$
\end{prop}\vspace{0.5cm}

This result follows from Theorem 3.3.4 of Nandagopalan (1990) or
from Ferreira (1994) by considering $T=1$ and replacing the
exceedance events $\{X_{i,j}>u_{n,j}\}$ by the upcrossing events
$\{X_{i,j}\leq u_{n,j}<X_{i+1,j}\}.$\vspace{0.5cm}

If ${\bf{S}}_n\xrightarrow [n\rain]{} {\bf{S}}[\alpha,\Pi]$ then,
for each $j=1,\ldots,d,$ $S_{n,j}\xrightarrow [n\rain]{} S_j$ where
$S_j\equiv S_j[\alpha_j,\Pi_j]$ is a compound Poisson process with
intensity $\alpha_j=\alpha(1-\Pi^*_j(0))$ and multiplicity
distribution
$$\Pi_j(\cdot)=\frac{\Pi_j^*(\cdot)}{1-\Pi_j^*(0)}\quad {\textrm{if}}\quad \Pi^*_j(0)<1,$$
where $\Pi_j^*$ represents the $j$-th projection of $\Pi.$ If
$\Pi_j^*(0)=1$ then $L_{S_j}(f)=1.$\vspace{0.5cm}

Suppose that for each ${\boldsymbol{\nu}}\in \R^d_+$ there exists
$\tilde{{\bf{u}}}^{({\boldsymbol{\nu}})}$ and that the corresponding
sequences of point processes
$\{{\bf{S}}^{({\bf{u}}^{(\boldsymbol{\nu})})}_n\}_{n\geq 1}$
converge to
${\bf{S}}[\alpha^{({\boldsymbol{\nu}})},\Pi^{({\boldsymbol{\nu}})}].$

Some questions naturally arise, at this point, about the way that
$\alpha$ and $\Pi$ depend on ${\boldsymbol{\nu}}.$

Suppose that the limit of
$\{{\bf{S}}^{({\bf{u}}^{(\boldsymbol{\nu})})}_n\}_{n\geq 1}$ doesn't
depend on the sequence $\tilde{{\bf{u}}}^{({\boldsymbol{\nu}})}$
considered and lets start by analyzing the effect of taking some
components of ${\boldsymbol{\nu}}$ converging to zero. From now on
we will generally write, for simplicity,
$\{{\bf{S}}^{({\boldsymbol{\nu}})}_n\}_{n\geq 1}$ instead of
$\{{\bf{S}}^{({\bf{u}}^{(\boldsymbol{\nu})})}_n\}_{n\geq
1}.$\vspace{0.3cm}

Let ${\boldsymbol{\nu}}'$ be such that $\nu_j=0$ if $j\in J\neq
\emptyset$ and $\nu_j>0$ for $j\in \{1,\ldots,d\}-J=\bar{J}.$ For
all ${\bf{y}}$ verifying $y_j>0$ for some $j\in J$ and all
$\tilde{{\bf{u}}}^{({\boldsymbol{\nu}'})},$
$$\lim_{n\to
\infty}P({\bf{S}}_n^{({\boldsymbol{\nu}}')}(A)={\bf{y}})=0$$ and,
for all ${\bf{y}}$ verifying $y_j=0$ for some $j\in J$ and any
$\tilde{{\bf{u}}}^{({\boldsymbol{\nu}'})},$
$$\lim_{n\to \infty}P\left(\bigcap_{j\in
\bar{J}}S_{n,j}^{(\nu_j')}(A)=y_j\right)-P({\bf{S}}_n^{({\boldsymbol{\nu}}')}(A)={\bf{y}})=0,\quad
\forall A\subset [0,1].$$\vspace{0.1cm}

\noindent Now, we have gathered the necessary conditions to prove
the following result.

\begin{prop}
Suppose that for each ${\boldsymbol{\nu}}\in \R_+^d$ there exist
$\tilde{{\bf{u}}}^{({\boldsymbol{\nu}})}$ and for each
$\tilde{{\bf{u}}}^{({\boldsymbol{\nu}})}$
${\bf{S}}_n^{({\boldsymbol{\nu}})}\xrightarrow
[n\rain]{}{\bf{S}}[\alpha^{({\boldsymbol{\nu}})},\Pi^{({\boldsymbol{\nu}})}].$
Let ${\boldsymbol{\nu}}'$  be such that $\nu_j=0$ if $j\in J\neq
\emptyset$ and $\nu_j>0$ if $j\in \bar{J}=\{1,\ldots,d\}-J.$ Define
$T({\boldsymbol{\nu}}')=\{{\boldsymbol{\nu}}\in \R^d_+:
\nu_j=\nu'_j\ {\textrm{for each } j\in \bar{J}}\}.$ Then
\begin{equation}
P({\bf{S}}[\alpha^{({\boldsymbol{\nu}}')},\Pi^{({\boldsymbol{\nu}}')}](A)=
{\bf{y}})=\lim_{\nu_j\to 0^+,\ j\in J \atop{{\boldsymbol{\nu}} \in
T({\boldsymbol{\nu}}')}}P({\bf{S}}[\alpha^{({\boldsymbol{\nu}})},\Pi^{({\boldsymbol{\nu}})}](A)={\bf{y}}),\label{eq2_5}
\end{equation}
for all $A\subset [0,1],$ ${\bf{y}}\in \N^d_0$ and
$$\Pi^{({\boldsymbol{\nu}}')}({\bf{y}})=\lim_{\nu_j\to 0^+,\ j\in J \atop{{\boldsymbol{\nu}} \in
T({\boldsymbol{\nu}}')}} \Pi^{({\boldsymbol{\nu}})}({\bf{y}}).$$
\end{prop}

\bdem Let $y_j=0$ for all $j\in J.$ Previously we have seen that
$$P({\bf{S}}[\alpha^{({\boldsymbol{\nu}}')},\Pi^{({\boldsymbol{\nu}}')}](A)=
{\bf{y}})=\lim_{n\to \infty}P\left(\bigcap_{j\in
\bar{J}}\{S_{n,j}^{(\nu_j')}(A)=y_j\}\right).$$ The following
inequalities also hold for all
$\tilde{{\bf{u}}}^{({\boldsymbol{\nu}})}$ with
${\boldsymbol{\nu}}\in T({\boldsymbol{\nu}}'),$
$$P\left(\bigcap_{j\in
\bar{J}}\{S_{n,j}^{(\nu_j')}(A)=y_j\}\right)\geq
P\left(\bigcap_{j\in \bar{J}}\{S_{n,j}^{(\nu_j')}(A)=y_j\},\
\bigcap_{j\in J} \{S_{n,j}^{(\nu_j)}(A)=0\}\right)$$ and
\begin{eqnarray*}
\lefteqn{\hspace{-2cm}-P\left(\bigcap_{j\in
\bar{J}}\{S_{n,j}^{(\nu_j')}(A)=y_j\},\ \bigcap_{j\in J}
\{S_{n,j}^{(\nu_j)}(A)=0\}\right)+P\left(\bigcap_{j\in
\bar{J}}\{S_{n,j}^{(\nu_j')}(A)=y_j\}\right)}\\[0.3cm]
&\leq& \sum_{j\in J}P(S_{n,j}^{(\nu_j)}(A)>0)\leq n \sum_{j\in
J}P(X_{1,j}\leq u_{n,j}^{(\nu_j)}<X_{2,j})\xrightarrow [n\rain]{}
\sum_{j\in J} \nu_j.
\end{eqnarray*}
Hence, $$\limsup_{n\to \infty}P\left(\bigcap_{j\in
\bar{J}}\{S_{n,j}^{(\nu_j')}(A)=y_j\}\right)\leq
P({\bf{S}}[\alpha^{({\boldsymbol{\nu}})},\Pi^{({\boldsymbol{\nu}})}](A)=
{\bf{y}})+\sum_{j\in J} \nu_j$$ and $$\liminf_{n\to
\infty}P\left(\bigcap_{j\in
\bar{J}}\{S_{n,j}^{(\nu_j')}(A)=y_j\}\right)\geq
P({\bf{S}}[\alpha^{({\boldsymbol{\nu}})},\Pi^{({\boldsymbol{\nu}})}](A)=
{\bf{y}}).$$ The result follows by taking limits when $\nu_j\to 0^+$
for all $j\in J.$\vspace{0.5cm}

\noindent Now, lets assume that $y_j>0$ for all $j\in J.$ Then,
$$\lim_{\nu_j\to 0^+,\ j\in J \atop{{\boldsymbol{\nu}} \in
T({\boldsymbol{\nu}}')}}P({\bf{S}}[\alpha^{({\boldsymbol{\nu}})},\Pi^{({\boldsymbol{\nu}})}](A)={\bf{y}})\leq
\lim_{\nu_j\to 0^+,\ j\in J \atop{{\boldsymbol{\nu}} \in
T({\boldsymbol{\nu}}')}}\lim_{n\to
\infty}P(S_{n,j}^{(\nu_j)}(A)>0)=0$$ and
$$P({\bf{S}}[\alpha^{({\boldsymbol{\nu}'})},\Pi^{({\boldsymbol{\nu}'})}](A)={\bf{y}})\leq
\lim_{n\to \infty}P(S_{n,j}^{(\nu_j)}(A)=y_j)=0,$$ so (\ref{eq2_5})
also holds in this case.

The second part of the result follows from the first one.\edem

The dependence  of the limit of ${\bf{S}}_n^{({\boldsymbol{\nu}})}$
on the sequence of levels disappears when we consider normalized
levels as in (\ref{eq2_2}) with the same number of exceedances
$\tau_j,$ $j=1,\ldots,d.$   This is, if ${\bf{u}}\equiv
\tilde{{\bf{u}}}^{({\boldsymbol{\nu}})}$ and ${\bf{v}}\equiv
\tilde{{\bf{v}}}^{({\boldsymbol{\nu}})}$ and furthermore
${\bf{u}}\equiv {\bf{u}}^{({\boldsymbol{\tau}})}$ and
${\bf{v}}\equiv {\bf{v}}^{({\boldsymbol{\tau}})}$ then the point
processes of upcrossings of ${\bf{u}}$ and of ${\bf{v}}$ have the
same limit since
$$P({\bf{S}}_n^{(\tilde{{\bf{u}}}^{({\boldsymbol{\nu}})})}(A)
\neq {\bf{S}}_n^{(\tilde{{\bf{v}}}^{({\boldsymbol{\nu}})})}(A))\leq
\sum_{j=1}^d|P(X_{i,j}>u_{n,j})-P(X_{i,j}>v_{n,j})|\xrightarrow
[n\rain]{} 0.$$

This last convergence holds if
$$nP(\min\{u_{n,j},v_{n,j}\}<X_{1,j}<\max \{u_{n,j},v_{n,j}\})\xrightarrow
[n\rain]{} 0,\quad j=1,\ldots,d.$$ Nevertheless, for some sequences
$\tilde{{\bf{u}}}^{({\boldsymbol{\nu}})},$ it's possible to relate
the limits of the corresponding point processes of upcrossings, as
proved in the following proposition.

\begin{prop}
Suppose ${\bf{X}}$ satisfies condition
$\Delta(\tilde{{\bf{u}}}^{({\boldsymbol{\nu}})})$ for all
$\tilde{{\bf{u}}}^{({\boldsymbol{\nu}})}.$ If for some
${\bf{u}}\equiv \tilde{{\bf{u}}}^{({\boldsymbol{\nu}})},$
${\bf{S}}_n^{({\boldsymbol{\nu}})}\xrightarrow
[n\rain]{}{\bf{S}}[\alpha^{({\boldsymbol{\nu}})},\Pi^{({\boldsymbol{\nu}})}]$
then, for any constant $c>0$ and
$\tilde{{\bf{u}}}^{(c{\boldsymbol{\nu}})}=\{{\bf{u}}^{({\boldsymbol{\nu}})}_{[n/c]}\}_{n\geq
1}$ we have
$${\bf{S}}_n^{(c{\boldsymbol{\nu}})}\xrightarrow
[n\rain]{}{\bf{S}}[c\alpha^{({\boldsymbol{\nu}})},\Pi^{({\boldsymbol{\nu}})}].$$
\end{prop}
\bdem We wish to prove that, for all $k>0,$
${\bf{y}}_1,\ldots,{\bf{y}}_k\in \N_0^d$ and $I_1,\ldots,I_k$
disjoint subintervals of $[0,1],$ we have
$$P\left(\bigcap_{j=1}^k\{{\bf{S}}_n^{(c{\boldsymbol{\nu}})}(I_j)={\bf{y}}_j\}\right)\xrightarrow
[n\rain]{}\prod_{j=1}^k
P({\bf{S}}[c\alpha^{({\boldsymbol{\nu}})},\Pi^{({\boldsymbol{\nu}})}](I_j)={\bf{y}}_j).$$
Since each interval $I_s$ can be written as the union of $m_s$
disjoint intervals of length less or equal than $\min\{1,c^{-1}\},$
from Lemma 2.1 we just have to prove that for an interval $I$ in
such conditions we have
$$P\left({\bf{S}}_n^{(c{\boldsymbol{\nu}})}(I)={\bf{y}}\right)\xrightarrow
[n\rain]{}P({\bf{S}}[c\alpha^{({\boldsymbol{\nu}})},\Pi^{({\boldsymbol{\nu}})}](I)={\bf{y}}),\quad
\forall {\bf{y}}\in \N_0^d.$$ Lets consider $I=(a,b],$ since for
another type of interval the arguments of the proof are analogous.

The stationarity assumption and ${\bf{u}}\equiv
\tilde{{\bf{u}}}^{(c{\boldsymbol{\nu}})}$ lead to $$\lim_{n\to
\infty}P({\bf{S}}_n^{(c{\boldsymbol{\nu}})}((a,b])={\bf{y}})=\lim_{n\to
\infty}P({\bf{S}}_n^{(c{\boldsymbol{\nu}})}((0,b-a])={\bf{y}}).$$

Comparing $\displaystyle{\lim_{n\to
\infty}P({\bf{S}}_n^{(c{\boldsymbol{\nu}})}((0,b-a])={\bf{y}})},$
for $\tilde{{\bf{u}}}^{(c{\boldsymbol{\nu}})},$ with
$\displaystyle{\lim_{n\to
\infty}P({\bf{S}}_{[n/c]}^{({\boldsymbol{\nu}})}((0,c(b-a)])={\bf{y}})},$
where
$$S_{[n/c],j}^{(\nu_j)}(\cdot)=\sum_{i=1}^{[n/c]}\indi_{\{X_{i,j}\leq u_{[n/c],j}^{(\nu_j)}<X_{i+1,j}\}}\delta_{\frac{i}{[n/c]}}(\cdot),\quad j=1,\ldots,d,$$
 we obtain
 \begin{eqnarray*}
 \lefteqn{\hspace{-1cm}\left|P\left(\bigcap_{j=1}^d \left\{\sum_{i=1}^{[nb-na]}\indi_{\{X_{i,j}\leq u_{[n/c],j}^{(\nu_j)}<X_{i+1,j}\}}=y_j
 \right\}\right)-P\left(\bigcap_{j=1}^d \left\{\sum_{i=1}^{\left[\left[\frac{n}{c}\right]c(b-a)\right]}\indi_{\{X_{i,j}
 \leq u_{[n/c],j}^{(\nu_j)}<X_{i+1,j}\}}=y_j \right\}
 \right)\right|}\\[0.3cm]
 &\leq& \sum_{j=1}^d P(X_{i,j}\leq
 u_{[n/c],j}^{(\nu_j)}<X_{i+1,j})=o(1).\hspace{7cm}
 \end{eqnarray*}
 Hence, we can conclude that
 \begin{eqnarray*}
 \lim_{n\to \infty}
 P({\bf{S}}_n^{(c{\boldsymbol{\nu}})}((a,b])={\bf{y}})&=&\lim_{n\to
 \infty}P({\bf{S}}_{[n/c]}^{({\boldsymbol{\nu}})}((ca,cb])={\bf{y}})\\[0.3cm]
 &=&
 P({\bf{S}}[\alpha^{({\boldsymbol{\nu}})},\Pi^{({\boldsymbol{\nu}})}]((ca,cb])={\bf{y}})\\[0.3cm]
&=&
P({\bf{S}}[c\alpha^{({\boldsymbol{\nu}})},\Pi^{({\boldsymbol{\nu}})}]((a,b])={\bf{y}}).
 \end{eqnarray*}
\edem

\section{The multivariate upcrossings index}

\pg Nandagopalan (1990) introduced a definition of a multivariate
extremal index for stationary sequences which we next recall. We
shall represent the vector of maxima of the first $n$ variables of
the i.i.d. sequence ${\hat{\bf{X}}},$ associated to ${\bf{X}},$ by
${\widehat{\bf{M}}}_n.$

\begin{defi}
A $d-$dimensional stationary sequence ${\bf X}$ is said to have
multivariate extremal index $\theta(\boldsymbol{\tau})$,
 $\boldsymbol{\tau} \in \R^d_+$, if
for each $\boldsymbol{\tau}=(\tau_1,\ldots,\tau_d)\in \R^d_+,$
exists ${\bf u}_n\equiv{\bf
u}_n^{(\boldsymbol{\tau})}=(u_{n1}^{(\tau_1)},\ldots,u_{nd}^{(\tau_d)}),\
n\geq 1,$ satisfying \vspace{0.3cm}

$$P(\widehat{\bf M}_n\leq {\bf u}_n)\xrightarrow
[n\rain]{}
 \widehat{\Psi}(\boldsymbol{\tau})
$$

\noindent and

$$P({\bf M}_n\leq {\bf u}_n)\xrightarrow
[n\rain]{}
(\widehat{\Psi}(\boldsymbol{\tau}))^{\theta(\boldsymbol{\tau})},\quad 0\leq
\theta(\boldsymbol{\tau})\leq 1
.$$
\end{defi}\vspace{0.5cm}

The existence of $0<\theta(\boldsymbol{\tau})<1$,
$\boldsymbol{\tau} \in
\R^d_+$, is an indicator of
the presence of clusters of events
$\displaystyle{\bigcup_{i=1}^d\{X_{i,j}>u_{n,j}\}}.$

In what follows we will extend the definition of upcrossings index
introduced in Ferreira (2006) to $d$-dimensional stationary
sequences,  relate this new index with  the multivariate extremal
index and present some of its properties.

\begin{defi}
${\bf{X}}$ has upcrossings index $\eta(\boldsymbol{\nu}),$
$\boldsymbol{\nu}\in \R^d_+,$ if for each $\boldsymbol{\nu}\in
\R^d_+$ exists $\tilde{\bf{u}}^{(\boldsymbol{\nu})},$
$$\lim_{n\to \infty}\exp\left(-nP\left(\bigcup_{j=1}^d\{X_{1,j}\leq
u_{n,j}^{(\nu_j)}<X_{2,j}\}\right)\right)=\varphi(\boldsymbol{\nu}),$$
and
$$\lim_{n\to \infty}P({\bf{S}}_n^{(\boldsymbol{\nu})}([0,1])={\bf{0}})=
(\varphi(\boldsymbol{\nu}))^{\eta(\boldsymbol{\nu})},\quad 0\leq
\eta(\boldsymbol{\nu})\leq 1,$$ for all
$\tilde{\bf{u}}^{(\boldsymbol{\nu})},$ $\boldsymbol{\nu}\in \R^d_+.$
\end{defi}

If for each $\boldsymbol{\nu}$ exist ${\bf{u}}\equiv
\tilde{{\bf{u}}}^{(\boldsymbol{\nu})}\equiv
{\bf{u}}^{(\boldsymbol{\tau})}$ then ${\bf{X}}$ has multivariate
extremal index $\theta(\boldsymbol{\tau})$ if and only if it has
multivariate upcrossings index $\eta(\boldsymbol{\nu})$ and, in this
case,
\begin{equation}
\eta(\boldsymbol{\nu})=\frac{\displaystyle{\lim_{n\to
\infty}nP\left(\bigcup_{j=1}^d
\{X_{1,j}>u_{n,j}\}\right)}}{\displaystyle{\lim_{n\to
\infty}nP\left(\bigcup_{j=1}^d \{X_{1,j}\leq
u_{n,j}<X_{2,j}\}\right)}}\ \theta(\boldsymbol{\tau}),\label{eq3_1}
\end{equation}
 since
$P(\widehat{\bf M}_n\leq {\bf u}_n)\xrightarrow [n\rain]{}
 \widehat{\Psi}(\boldsymbol{\tau})$ if and only if $nP({\bf{X}}_1\nleqslant {\bf u}_n)\xrightarrow
[n\rain]{} -\log \widehat{\Psi}(\boldsymbol{\tau}).$\vspace{0.5cm}

The relation (\ref{eq3_1}) generalizes the one obtained by Ferreira
(2006) for $d=1$ and, as we will see next, it allows us to obtain
similar relations for the marginal upcrossings indices and extremal
indices.

The upcrossings index of ${\bf{X}}_j,$ $j=1,\ldots,d,$ can be easily
obtained from the knowledge of the multivariate upcrossings index,
as stated in the next result. Nevertheless, for some sequences the
multivariate upcrossings index can also be obtained from the
marginal upcrossings indices, as we will see in the next section.

\begin{prop}
If ${\bf{X}}$ has upcrossings index
$\eta^{{\bf{X}}}(\boldsymbol{\nu})>0,$ $\boldsymbol{\nu} \in
\R^d_+,$ then, for any $\emptyset\neq J=\{i_1,\ldots,i_s\}\subset
\{1,\ldots,d\},$
${\bf{X}}_{i_1,\ldots,i_s}=\{(X_{n,i_1},\ldots,X_{n,i_s})\}_{n\geq
1}$ has upcrossings index
$$\eta^{{\bf{X}}_{i_1,\ldots,i_s}}(\boldsymbol{\nu}')=\lim_{\nu_j\to 0^+,\ j\in \bar{J}
\atop{\boldsymbol{\nu}\in
\overline{T}(\boldsymbol{\nu}')}}\eta^{{\bf{X}}}(\boldsymbol{\nu}),\quad
\boldsymbol{\nu}'\in \R^s_+,$$ with
$\overline{T}(\boldsymbol{\nu}')=\{\boldsymbol{\nu} \in \R^d_+:\
\nu_j=\nu_j',\ j\in J\}$ and $\overline{J}=\{1,\ldots,d\}-J.$
\end{prop}
\bdem For each $\boldsymbol{\nu}' \in \R^s_+$ consider
$\boldsymbol{\nu}'' \in \R^d_+$ such that
$$\nu_j''=\nu_j'\quad {\textrm{if}}\ j\in J\ {\textrm{and}}\
\nu_j''=0\ {\textrm{if}}\ j\in \bar{J}.$$ We have
$\overline{T}(\boldsymbol{\nu}')=\overline{T}(\boldsymbol{\nu}'')$
and, as previously seen, by considering levels
$\tilde{{\bf{u}}}^{(\boldsymbol{\nu}'')}$ and
$\tilde{{\bf{u}}}^{(\boldsymbol{\nu})},$ $\boldsymbol{\nu}\in
\overline{T}(\boldsymbol{\nu}''),$
\begin{eqnarray*}
\lim_{n\to
\infty}P\left(\bigcap_{j=1}^s\{S_{n,i_j}^{(\nu'_{i_j})}([0,1])=0\}\right)&=&\lim_{n\to
\infty}P\left(\bigcap_{j=1}^d\{
S_{n,j}^{(\nu''_j)}([0,1])=0\}\right)\\[0.3cm]
&=&\lim_{\nu_j\to 0^+,\ j\in \bar{J} \atop{\boldsymbol{\nu}\in
\overline{T}(\boldsymbol{\nu}'')}} \lim_{n\to \infty}
P({\bf{S}}_n^{(\boldsymbol{\nu})}([0,1])={\bf{0}})
\end{eqnarray*}
and
\begin{eqnarray*}
\lim_{n\to \infty} nP\left(\bigcup_{j=1}^s\{ X_{1,i_j}\leq
u_{n,i_j}^{(\nu'_{i_j})}<X_{2,i_j}\}\right)&=&\lim_{n\to \infty}
nP\left(\bigcup_{j=1}^d\{ X_{1,j}\leq
u_{n,j}^{(\nu''_j)}<X_{2,j}\}\right)\\[0.3cm]
&=& \lim_{\nu_j\to 0^+,\ j\in \bar{J} \atop{\boldsymbol{\nu}\in
\overline{T}(\boldsymbol{\nu}'')}} \lim_{n\to \infty}
nP\left(\bigcup_{j=1}^d\{ X_{1,j}\leq
u_{n,j}^{(\nu_j)}<X_{2,j}\}\right).
\end{eqnarray*}

Hence, for ${\bf{X}}_{i_1,\ldots,i_s}$ there exist levels
$\tilde{{\bf{u}}}^{(\boldsymbol{\nu}')},$ $\boldsymbol{\nu}'\in
\R^s_+,$ such that $$\lim_{n\to \infty}
\exp\left(-nP\left(\bigcup_{j=1}^s \{X_{1,i_j}\leq
u_{n,i_j}^{(\nu'_{i_j})}<X_{2,i_j}\}\right)\right)=\xi(\boldsymbol{\nu}')$$
and $$\lim_{n\to
\infty}P({\bf{S}}_n^{(\boldsymbol{\nu}')}([0,1])={\bf{0}})=(\xi(\boldsymbol{\nu}'))^{\displaystyle{\lim_{\nu_j\to
0^+,\ j\in \bar{J} \atop{\boldsymbol{\nu}\in
\overline{T}(\boldsymbol{\nu}')}}\eta(\boldsymbol{\nu})}}.$$\edem

We will say that ${\bf{X}}$ has an upcrossing at $i$ when the event
$\displaystyle{\bigcup_{j=1}^d\{X_{i,j}\leq u_{n,j}<X_{i+1,j}\}}$
occurs.

The next result states that the reciprocal of the multivariate
upcrossings index can be interpreted as the limiting mean cluster
size of upcrossings of ${\bf{X}},$ which generalizes Proposition 2
in Ferreira (2007).

\begin{prop}
Suppose that for each $\boldsymbol{\nu}\in \R^d_+$ there exist
levels $\tilde{{\bf{u}}}^{(\boldsymbol{\nu})}$ and ${\bf{X}}$
satisfies condition $\Delta(\tilde{{\bf{u}}}^{(\boldsymbol{\nu})}).$
If ${\bf{X}}$ has multivariate upcrossings index
$\eta(\boldsymbol{\nu})>0,$ $\boldsymbol{\nu}\in \R^d_+,$ then
$$\sum_{k\geq 1}k\overline{\Pi}^{(\boldsymbol{\nu})}(k)\xrightarrow
[n\rain]{}\frac{1}{\eta(\boldsymbol{\nu})},\quad  {\textrm{for
all}}\  \tilde{{\bf{u}}}^{(\boldsymbol{\nu})},\ \boldsymbol{\nu}\in
\R^d_+,$$ where
$$\overline{\Pi}^{(\boldsymbol{\nu})}(k)=P\left(\sum_{i=1}^{r_n}
\indi_{\bigcup_{j=1}^d\{ X_{i,j}\leq
u_{n,j}^{(\nu_j)}<X_{i+1,j}\}}=k\ \Bigg| \sum_{i=1}^{r_n}
\indi_{\bigcup_{j=1}^d \{X_{i,j}\leq
u_{n,j}^{(\nu_j)}<X_{i+1,j}\}}>0\right).$$
\end{prop}
\bdem  Noting that
$$P({\bf{S}}_n^{(\boldsymbol{\nu})}(C_{n,1})={\bf 0})=P\left(\sum_{i=1}^{r_n}
\indi_{\bigcup_{j=1}^d \{X_{i,j}\leq
u_{n,j}^{(\nu_j)}<X_{i+1,j}\}}=0\right),$$ and, from Lemma 2.1,
\begin{eqnarray*}
\lefteqn{\hspace{-3cm}\lim_{n\to \infty}k_nP\left(\sum_{i=1}^{r_n}
\indi_{\bigcup_{j=1}^d\{ X_{i,j}\leq
u_{n,j}^{(\nu_j)}<X_{i+1,j}\}}>0\right)}\\[0.3cm]
&=&\lim_{n\to
\infty}k_nP({\bf{S}}_n^{(\boldsymbol{\nu})}(C_{n,1})\neq{\bf 0})\\[0.3cm]
&=& -\log \left(\lim_{n\to
\infty}P({\bf{S}}_n^{(\boldsymbol{\nu})}([0,1])={\bf 0})\right)\\[0.3cm]
&=& \eta(\boldsymbol{\nu})\log
\varphi(\boldsymbol{\nu}),\hspace{2cm}
\end{eqnarray*}
hence
\begin{eqnarray*}
\sum_{k\geq
1}k\overline{\Pi}^{(\boldsymbol{\nu})}(k)&=&\frac{\displaystyle{\left[\frac{n}{k_n}\right]P\left(\bigcup_{j=1}^d
\{X_{1,j}\leq
u_{n,j}^{(\nu_j)}<X_{2,j}\}\right)}}{\displaystyle{P\left(\sum_{i=1}^{r_n}
\indi_{\bigcup_{j=1}^d \{X_{i,j}\leq
u_{n,j}^{(\nu_j)}<X_{i+1,j}\}}>0\right)}}\\[0.3cm]
&\xrightarrow [n\rain]{}& \frac{\log
\varphi(\boldsymbol{\nu})}{\eta(\boldsymbol{\nu})\log
\varphi(\boldsymbol{\nu})}=\frac{1}{\eta(\boldsymbol{\nu})}.
\end{eqnarray*}\edem

\section{Computation of the multivariate upcrossings index from the uni\-va\-riate indices}

\pg As shown in the previous section, if ${\bf{X}}$ has multivariate
upcrossings index $\eta(\boldsymbol{\nu}),$ $\boldsymbol{\nu}\in
\R^d_+,$ then each marginal sequence ${\bf{X}}_j,$ $j=1,\ldots,d,$
has upcrossings index $\displaystyle{\eta_j=\lim_{\nu_i\to 0^+,\
i\neq j}\eta(\boldsymbol{\nu})}.$ Since $1-$dimensional sequences
are more easy to handle it would be interesting to obtain the
multivariate index from the marginal indices $\eta_j.$ This can be
done if the sequence ${\bf{X}}$ satisfies asymptotic independence
conditions on its marginals. Several of conditions of this type have
been introduced in the literature during the last years (Davis
(1982), H\"usler (1990), Ferreira (1994), Pereira (2002)),
simplifying the study of multivariate extremes through the study of
the univariate ones.

In this section we introduce a new condition of asymptotic
independence that allows us to handle componentwise the upcrossings of
${\bf{X}}$ and hence simplify the computation of
$\eta(\boldsymbol{\nu}).$ In the next section we shall see how this
type of condition can be weakened in the presence of other
conditions that restrict rapid oscillations of
${\bf{X}}_j.$\vspace{0.5cm}

\begin{defi}
${\bf{X}}$ satisfies condition $H({\bf{u}})$ when
$$\sum_{1\leq j<j'\leq d} n\sum_{i=1}^{[n/k_n]} P(X_{1,j}\leq u_{n,j}<X_{2,j},
\ X_{i,j'}\leq u_{n,j'}<X_{i+1,j'})\xrightarrow [n\rain]{}0,$$ for
some sequence ${\bf{k}}=\{k_n\}_{n\geq 1}$ such that
$k_n\xrightarrow [n\rain]{}0$ and $k_n/n\xrightarrow [n\rain]{}0.$
\end{defi}\vspace{0.3cm}

\begin{prop}
Suppose that for each $\boldsymbol{\nu}\in \R^d_+$ there exist
levels $\tilde{{\bf{u}}}^{(\boldsymbol{\nu})}$  and ${\bf{X}}$
satisfies condition $\Delta(\tilde{{\bf{u}}}^{(\boldsymbol{\nu})})$
and condition $H(\tilde{{\bf{u}}}^{(\boldsymbol{\nu})}),$ for all
$\tilde{{\bf{u}}}^{(\boldsymbol{\nu})}$ and some sequence ${\bf{k}}$
satisfying (\ref{eq2_4}).

If ${\bf{X}}_j$ has upcrossings index $\eta_j,$ $j=1,\ldots,d,$ then
\begin{enumerate}
\item[{\bf{(i)}}] ${\bf{X}}$ has multivariate upcrossings index $$\eta(\boldsymbol{\nu})=
\frac{\displaystyle{\sum_{i=1}^d
\eta_j\nu_j}}{\displaystyle{\sum_{i=1}^d \nu_j}},\quad
\boldsymbol{\nu}\in \R^d_+.
$$
\item[{\bf{(ii)}}] For all ${\bf{y}}$ with two or more non zero components
and all $\tilde{{\bf{u}}}^{(\boldsymbol{\nu})},$
$\Pi_n^{(\boldsymbol{\nu})}({\bf{y}})\xrightarrow [n\rain]{}0.$
\item[{\bf{(iii)}}] For all $y>0$ and
$\tilde{{\bf{u}}}^{(\boldsymbol{\nu})},$
$$\Pi_n^{(\boldsymbol{\nu})}(0,\ldots,0,y,0,\ldots,0)-\frac{\nu_j\eta_j}
{\sum_{j=1}^d\nu_j\eta_j}\ P(S_n^{(\nu_j)}(C_{n,1})=y\ |\
S_n^{(\nu_j)}(C_{n,1})>0)\xrightarrow [n\rain]{}0.$$
\end{enumerate}
\end{prop}
\bdem {\bf{(i)}} Under condition
$\Delta(\tilde{{\bf{u}}}^{(\boldsymbol{\nu})})$ we have
$$\lim_{n\to \infty} P({\bf{S}}_n^{(\boldsymbol{\nu})}([0,1])={\bf{0}})=
\lim_{n\to \infty} \exp \left(-k_n
P\left(\sum_{i=1}^{r_n}\indi_{\bigcup_{j=1}^d \{X_{i,j}\leq
u_{n,j}^{(\nu_j)}<X_{i+1,j}\}}>0\right)\right),$$ being this last
limit equal to
$$\lim_{n\to \infty}\exp \left(-k_n\sum_{j=1}^d P\left(\sum_{i=1}^{r_n}\indi_{\{X_{i,j}\leq
u_{n,j}^{(\nu_j)}<X_{i+1,j}\}}>0\right)\right)
$$ since $$P\left(\sum_{i=1}^{r_n}\indi_{\bigcup_{j=1}^d \{X_{i,j}\leq
u_{n,j}^{(\nu_j)}<X_{i+1,j}\}}>0\right)\leq
P\left(\bigcup_{j=1}^d\left\{\sum_{i=1}^{r_n}\indi_{\{X_{i,j}\leq
u_{n,j}^{(\nu_j)}<X_{i+1,j}\}}>0\right\}\right)$$ and
\begin{eqnarray*}
P\left(\sum_{i=1}^{r_n}\indi_{\bigcup_{j=1}^d \{X_{i,j}\leq
u_{n,j}^{(\nu_j)}<X_{i+1,j}\}}>0\right)\geq
P\left(\bigcup_{j=1}^d\left\{\sum_{i=1}^{r_n}\indi_{\{X_{i,j}\leq
u_{n,j}^{(\nu_j)}<X_{i+1,j}\}}>0\right\}\right)\hspace{3cm}\\[0.3cm]
-P\left(\bigcup_{1\leq j<j'\leq d}\left\{\sum_{1\leq i,i'\leq
r_n}\indi_{\{X_{i,j}\leq u_{n,j}^{(\nu_j)}<X_{i+1,j},\ X_{i',j'}\leq
u_{n,j'}^{(\nu_{j'})}<X_{i'+1,j'}\}}>0\right\}\right),
\end{eqnarray*}
where the last term is bounded by the sum that defines condition
$H(\tilde{{\bf{u}}}^{(\boldsymbol{\nu})}).$

Then, for all $\tilde{{\bf{u}}}^{(\boldsymbol{\nu})},$ we have
\begin{eqnarray*}
\lim_{n\to \infty}
P({\bf{S}}_n^{(\boldsymbol{\nu})}([0,1])={\bf{0}})&=&
\exp\left(\sum_{j=1}^d \log\left(\lim_{n\to
\infty}P(S_n^{(\nu_j)}([0,1])=0)\right)\right)\\[0.3cm]
&=& \exp\left(-\sum_{j=1}^d \eta_j \nu_j\right)
\end{eqnarray*}
and $\varphi(\boldsymbol{\nu})=\exp\left(-\sum_{j=1}^d
\nu_j\right),$ since from the asymptotic behaviour of the first term
of $H(\tilde{{\bf{u}}}^{(\boldsymbol{\nu})}),$ $$\lim_{n\to \infty}n
P\left(\bigcup_{j=1}^d \{X_{1,j}\leq
u_{n,j}^{(\nu_j)}<X_{2,j}\}\right)=\lim_{n\to \infty} \sum_{j=1}^d
nP(X_{1,j}\leq u_{n,j}^{(\nu_j)}<X_{2,j})=\sum_{j=1}^d
\nu_j.$$\vspace{0.3cm}

\noindent {\bf{(ii)}} and {\bf{(iii)}} follow by applying
$H(\tilde{{\bf{u}}}^{(\boldsymbol{\nu})}),$ in an analogous way to
what has been done previously and by noting again that
$k_nP(S_n^{(\nu_j)}(I_{n,1})>0)\xrightarrow [n\rain]{} \nu_j\eta_j,$
$j=1,\ldots,d.$ \edem

\section{Mild oscillation restriction}

\pg If we consider on $A_{i,n}=\bigcup_{j=1}^d\{X_{i,j}\leq
u_{n,j}<X_{i+1,j}\}$, $i\geq 1$,  analogous restrictions to the ones considered
in Ferreira (2006) and denoted by $\tilde{D}^{(k)}({\bf{u}}),$ we
find formulas to compute $\eta(\boldsymbol{\nu})$ from a finite
number of variables of ${\bf{X}}.$ Since no new arguments are needed
to obtain such results we shall restrict ourselves to
$\tilde{D}^{(3)}({\bf{u}})$ which will be applied in the following
examples.

The  local dependence condition $\tilde{D}^{(3)}({\bf{u}})$ also
allows us to consider the condition  $H({\bf{u}})$ with a finite number
of terms.\vspace{0.5cm}

\begin{prop}
Suppose that for each $\boldsymbol{\nu}\in \R^d_+$  there exist
levels $ \tilde{{\bf{u}}}^{(\boldsymbol{\nu})}$ and ${\bf{X}}$
satisfies condition $\Delta(\tilde{{\bf{u}}}^{(\boldsymbol{\nu})})$
and
\begin{equation}
nP\left(A_{1,n}^{(\boldsymbol{\nu})},\
\overline{A}_{2,n}^{(\boldsymbol{\nu})},\
\overline{A}_{3,n}^{(\boldsymbol{\nu})},\
\sum_{i=4}^{[n/k_n]}\indi_{A_{i,n}^{(\boldsymbol{\nu})}}>0\right)\xrightarrow
[n\rain]{}0,\label{eq5_1}
\end{equation}
for all $ \tilde{{\bf{u}}}^{(\boldsymbol{\nu})}$ and some sequence ${\bf{k}}$ satisfying (\ref{eq2_4}). Then
${\bf{X}}$ has upcrossings index $\eta(\boldsymbol{\nu}),$
$\boldsymbol{\nu}\in \R^d_+,$ if and only if
$$\frac{P(A_{1,n}^{(\boldsymbol{\nu})},\ \overline{A}_{2,n}^{(\boldsymbol{\nu})},\
\overline{A}_{3,n}^{(\boldsymbol{\nu})})}{P(A_{1,n}^{(\boldsymbol{\nu})})}\xrightarrow
[n\rain]{}\eta(\boldsymbol{\nu}),$$ for all
$\tilde{{\bf{u}}}^{(\boldsymbol{\nu})},$ with
$\displaystyle{A_{i,n}^{(\boldsymbol{\nu})}=\bigcup_{j=1}^d\{X_{i,j}\leq
u_{n,j}^{(\nu_j)}<X_{i+1,j}\}},$ $i\geq 1.$
\end{prop}

This result is proved with analogous arguments to the ones used in
 Proposition 3.1 and Corollary 3.1 of
Ferreira (2006) and therefore its proof  will be
omitted.\vspace{0.3cm}

As a consequence of the two previous propositions we can now obtain
new ways to compute $\eta(\boldsymbol{\nu}).$

\begin{coro}
Suppose that for each $\boldsymbol{\nu}\in \R^d_+$ there exist
levels $ \tilde{{\bf{u}}}^{(\boldsymbol{\nu})},$  ${\bf{X}}$
satisfies the conditions
$\Delta(\tilde{{\bf{u}}}^{(\boldsymbol{\nu})})$, (\ref{eq5_1}) and
$$\sum_{1\leq j<j'\leq d} n \sum_{i=1}^3P(X_{1,j}\leq u_{n,j}^{(\nu_j)}<X_{2,j},
\ X_{i,j'}\leq u_{n,j'}^{(\nu_{j'})}<X_{i+1,j'})\xrightarrow
[n\rain]{} 0,$$
 for all $ \tilde{{\bf{u}}}^{(\boldsymbol{\nu})}$
and some sequence ${\bf{k}}$ satisfying
(\ref{eq2_4}). Then ${\bf{X}}$ has upcrossings index
$\eta(\boldsymbol{\nu}),$ $\boldsymbol{\nu}\in \R^d_+,$ if and only
if ${\bf{X}}_j$ has upcrossings index $\eta_j,$ $j=1,\ldots,d,$ and
in this case
$$\eta(\boldsymbol{\nu})=\frac{\displaystyle{\sum_{j=1}^d \nu_j\eta_j}}{\displaystyle{\sum_{j=1}^d \nu_j}},$$
with $\eta_j=\lim_{n\to \infty}P(A_{1,n}^{(\nu_j)},\
\overline{A}_{3,n}^{(\nu_j)})/P(A_{1,n}^{(\nu_j)}),$ $j=1,\ldots,d,$
$\displaystyle{A_{i,n}^{(\nu_j)}=\{X_{i,j}\leq
u_{n,j}^{(\nu_j)}<X_{i+1,j}\}}.$
\end{coro}

\begin{coro}
Suppose that for each $\boldsymbol{\nu}\in \R^d_+$ there exist
levels $ \tilde{{\bf{u}}}^{(\boldsymbol{\nu})},$ ${\bf{X}}$ verifies
 $\Delta(\tilde{{\bf{u}}}^{(\boldsymbol{\nu})}),$
$H(\tilde{{\bf{u}}}^{(\boldsymbol{\nu})})$ and
$$nP\left(A_{1,n}^{(\nu_j)},\ \overline{A}_{3,n}^{(\nu_j)},\ \sum_{i=4}^{[n/k_n]}
\indi_{A_{i,n}^{(\nu_j)}}>0\right)\xrightarrow [n\rain]{} 0,\quad
j=1,\ldots,d,$$ for all $ \tilde{{\bf{u}}}^{(\boldsymbol{\nu})}$ and some sequence ${\bf{k}}$ satisfying
(\ref{eq2_4}).

Then ${\bf{X}}$ has multivariate upcrossings index
$\eta(\boldsymbol{\nu}),$ $\boldsymbol{\nu}\in \R^d_+,$ if and only
if $$\lim_{n\to \infty}\frac{P(A_{1,n}^{(\nu_j)},\
\overline{A}_{3,n}^{(\nu_j)})}{P(A_{1,n}^{(\nu_j)})}\xrightarrow
[n\rain]{}\eta_j,\quad j=1,\ldots,d,$$ and, in this case,
$$\eta(\boldsymbol{\nu})=\frac{\sum_{j=1}^d\nu_j\eta_j}{\sum_{j=1}^d\nu_j},\quad \boldsymbol{\nu}\in \R^d_+.$$
\end{coro}

\section{Examples}

\begin{ex}
\rm{Let ${\bf{Y}}=\{Y_n\}_{n\geq -2}$ be a sequence of independent
random  variables  uniformly distributed on [0,1], $u_n=1-\tau'/n,$
$\tau'>0, \ n\geq 1.$

Consider the 2-dimensional sequence
${\bf{X}}=\{(X_{n,1},X_{n,2})\}_{n\geq 1}$ with $X_{n,1}=\max\{Y_n,
Y_{n-2}, Y_{n-3}\}$ and $X_{n,2}=Y_{n+1},$ $n\geq 1.$ For
${\bf{u}}=\{(u_{n,1},u_{n,2})\}_{n\geq 1}$ with $u_{n,1}=1-\tau'_1/n$
and $u_{n,2}=1-\tau'_2/n,$ $n\geq 1,$ we have ${\bf{u}}\equiv
{\bf{u}}^{(\boldsymbol{\tau})}$ with
$\boldsymbol{\tau}=(3\tau'_1,\tau'_2)$ and  ${\bf{u}}\equiv
\tilde{{\bf{u}}}^{(\boldsymbol{\nu})}$ with
$\boldsymbol{\nu}=(2\tau'_1,\tau'_2)$.

${\bf{X}}$ satisfies condition $\Delta({\bf{u}})$ since it's
4-dependent and condition $H({\bf{u}})$ since
$$nP(X_{1,1}\leq u_{n,1}<X_{2,1},\ X_{1,2}\leq u_{n,2}<X_{2,2})\leq nP(Y_0>u_{n,1},\ Y_3>u_{n,2})\xrightarrow
[n\rain]{} 0,$$ and, for $i\geq 2$ and ${\bf{k}}$ verifying
(\ref{eq2_4}),
\begin{eqnarray*}
\lefteqn{\hspace{-2cm}\left[\frac{n}{k_n}\right]nP(X_{1,1}\leq
u_{n,1}<X_{2,1},\
X_{i,1}\leq u_{n,2}<X_{i+1,2})}\\[0.3cm]
&\leq& \left[\frac{n}{k_n}\right]nP(Y_0>u_{n,1}\ \vee Y_2>u_{n,1},\
Y_{i+1}>u_{n,2})\xrightarrow [n\rain]{} 0.
\end{eqnarray*}
Hence, from the study in Sebastiao (2007) for the first marginal we
conclude that ${\bf{X}}_1$ has upcrossings index equal to 1/2,
${\bf{X}}_2$ has upcrossings index equal to 1 and ${\bf{X}}$ has
multivariate upcrossings index
$$\eta(\boldsymbol{\nu})=\frac{\nu_1/2+\nu_2}{\nu_1+\nu_2},\quad \boldsymbol{\nu}\in \R^2_+,$$
and multivariate extremal index
$$\theta(\boldsymbol{\tau})=\frac{\nu_1+\nu_2}{\tau_1+\tau_2}\ \eta(\boldsymbol{\nu})=
\frac{\tau_1/3+\tau_2}{\tau_1+\tau_2}, \quad \boldsymbol{\tau}\in
\R^2_+.$$ Furthermore,
\begin{flushleft}
\begin{tabular}{l}
$\lim_{n\to \infty} \Pi_n^{(\nu_1)}(2)=\lim_{n\to
\infty}\Pi_n^{(\nu_2)}(1)=1,\qquad \lim_{n\to
\infty}\Pi_n^{(\nu_1,\nu_2)}(2,0)=\frac{\nu_1/2}{\nu_1/2+\nu_2}$
\\[0.2cm] $\lim_{n\to
\infty}\Pi_n^{(\nu_1,\nu_2)}(0,1)=\frac{\nu_2}{\nu_1/2+\nu_2}.$
\end{tabular}
\end{flushleft}}
\end{ex}\vspace{0.5cm}

\begin{ex}
\rm{ Lets consider now the 2-dimensional sequence
${{\bf{Z}}=\{Z_{n,1},Z_{n,2}\}_{n\geq 1}}$ with $Z_{n,1}=X_{n,1}$
and $Z_{n,2}=Y_n,$ $n\geq 1$. For the same sequence of levels
${\bf{u}}\equiv{\bf{u}}^{(\boldsymbol{\tau})}\equiv
\tilde{{\bf{u}}}^{(\boldsymbol{\nu})},$ condition  $H({\bf{u}})$
doesn't hold since $$\lim_{n\to \infty}nP(Z_{1,1}\leq
u_{n,1}<Z_{2,1},\ Z_{1,2}\leq u_{n,2}<Z_{2,2})=\lim_{n\to
\infty}nP(Y_2>\max\{u_{n,1},u_{n,2}\})=\min_{j=1,2}\tau_j>0.$$

\noindent ${\bf{Z}}$ verifies condition (\ref{eq5_1}) since each of
its marginals verifies $\tilde{D}^{(3)}({\bf{u}}_j)$ and
$$nP\left(A_{1,n}^{(\nu_j)}, \overline{A}_{3,n}^{(\nu_j)},\ \bigcap_{i=1}^3 \overline{A}_{i,n}^{(\nu_{j'})},
\ \sum_{i=4}^{r_n}\indi_{A_{i,n}^{(\nu_{j'})}}>0\right)\leq 2n
P(Y_1>u_{n,j})P(Y_1>u_{n,j'})=o(1).$$

\noindent We can now compute the multivariate upcrossings index from
Proposition 5.1.
\begin{eqnarray*}
nP\left(\bigcup_{j=1}^d\{Z_{1,j}\leq
u_{n,j}<Z_{2,j}\}\right)=nP(A_{1,n}^{(\nu_1)},A_{1,n}^{(\nu_2)})+nP(A_{1,n}^{(\nu_1)},\overline{A}_{1,n}^{(\nu_2)})+
nP(\overline{A}_{1,n}^{(\nu_1)},A_{1,n}^{(\nu_2)}),
\end{eqnarray*}
where
\begin{eqnarray*}
\lim_{n\to
\infty}nP(A_{1,n}^{(\nu_1)},A_{1,n}^{(\nu_2)})&=&\lim_{n\to
\infty}nP(Y_2>\max\{u_{n,1},u_{n,2}\})= \min_{j=1,2}\tau_j,\\[0.3cm]
\lim_{n\to
\infty}nP(A_{1,n}^{(\nu_1)},\overline{A}_{1,n}^{(\nu_2)})&=&\lim_{n\to
\infty}nP(Y_0>u_{n,1})+nP(u_{n,1}<Y_2\leq u_{n,2}) =\left\{
\begin{array}{lcl}
\tau_1 & \textrm{if} & \tau_2\geq \tau_1\\
2\tau_1-\tau_2 & \textrm{if} & \tau_2< \tau_1,
\end{array}\right.\\[0.3cm]
\lim_{n\to
\infty}nP(\overline{A}_{1,n}^{(\nu_1)},A_{1,n}^{(\nu_2)})&=&\left\{
\begin{array}{lcl}
\tau_2-\tau_1 & \textrm{if} & \tau_2\geq \tau_1\\
0 & \textrm{if} & \tau_2< \tau_1.
\end{array}\right.
\end{eqnarray*}

\noindent Hence, $$\varphi({\boldsymbol{\nu}})=\left\{
\begin{array}{lcl}
e^{-\left(\frac{\nu_1}{2}+\nu_2\right)} & \textrm{if} & 2\nu_2\geq \nu_1\\
e^{-\nu_1} & \textrm{if} & 2\nu_2< \nu_1.
\end{array}\right.$$

\noindent Moreover,
\begin{eqnarray*}
&&\lim_{n\to\infty}nP(A_{1,n}^{(\nu_1)},A_{1,n}^{(\nu_2)},\overline{A}_{3,n}^{(\boldsymbol{\nu})})=0,\\[0.2cm]
&&\lim_{n\to\infty}nP(A_{1,n}^{(\nu_1)},\overline{A}_{1,n}^{(\nu_2)},\overline{A}_{2,n}^{(\nu_2)},
\overline{A}_{3,n}^{(\boldsymbol{\nu})})=\lim_{n\to
\infty}nP(Y_0>u_{n,1})=\tau_1,\\[0.2cm]
&&\lim_{n\to\infty}nP(\overline{A}_{1,n}^{(\nu_1)},A_{1,n}^{(\nu_2)},\overline{A}_{2,n}^{(\nu_1)},
\overline{A}_{3,n}^{(\boldsymbol{\nu})})=\lim_{n\to
\infty}nP(u_{n,2}<Y_2\leq_{u_n,1})=\left\{
\begin{array}{lcl}
\tau_2-\tau_1 & \textrm{if} & \tau_2\geq \tau_1\\
0 & \textrm{if} & \tau_2< \tau_1.
\end{array}\right.
\end{eqnarray*}

\noindent We then obtain, $$\eta(\boldsymbol{\nu})=\lim_{n\to
\infty}\frac{nP(A_{1,n}^{(\boldsymbol{\nu})},\overline{A}_{2,n}^{(\boldsymbol{\nu})},
\overline{A}_{3,n}^{(\boldsymbol{\nu})})}{nP(A_{1,n}^{(\boldsymbol{\nu})})}=\left\{
\begin{array}{lcl}
\frac{\nu_2}{\nu_1/2+\nu_2} & \textrm{if} & 2\nu_2\geq \nu_1\\
1/2 & \textrm{if} & 2\nu_2< \nu_1,
\end{array}\right.$$
and consequently $\eta_1=1/2$ and $\eta_2=1$, as expected.

As stated in Proposition 3.2, we have
\begin{eqnarray*}
\sum_{k=1}^2kP\left(\sum_{i=1}^{r_n} \indi_{\bigcup_{j=1}^2\{
X_{i,j}\leq u_{n,j}^{(\nu_j)}<X_{i+1,j}\}}=k\ \Bigg|
\sum_{i=1}^{r_n} \indi_{\bigcup_{j=1}^2 \{X_{i,j}\leq
u_{n,j}^{(\nu_j)}<X_{i+1,j}\}}>0\right)\\[0.3cm]
=1\times
\Pi_n^{(\nu_1,\nu_2)}(0,1)+2(\Pi_n^{(\nu_1,\nu_2)}(2,0)+\Pi_n^{(\nu_1,\nu_2)}(2,1))\xrightarrow
[n\rain]{} \frac{1}{\eta(\boldsymbol{\nu})},
\end{eqnarray*}

\noindent since
\begin{eqnarray*} \lim_{n\to
\infty}\Pi_n^{(\nu_1,\nu_2)}(2,0)&=&\lim_{n\to
\infty}\frac{nP(u_{n,1}<Y_2\leq
u_{n,2})}{k_nP({\bf{S}}_n^{(\boldsymbol{\nu})}(C_{n,1})\neq
{\bf{0}})}=\left\{
\begin{array}{lcl}
0 & \textrm{if} & \nu_1\leq 2\nu_2\\
\frac{\nu_1/2-\nu_2}{\nu_1/2} & \textrm{if} & \nu_1> 2\nu_2,
\end{array}\right.\\[0.3cm]
\lim_{n\to \infty}\Pi_n^{(\nu_1,\nu_2)}(0,1)&=&\lim_{n\to
\infty}\frac{nP(u_{n,2}<Y_2\leq u_{n,1})}{-\log
\varphi(\boldsymbol{\nu})\eta(\boldsymbol{\nu})}=\left\{
\begin{array}{lcl}
\frac{\nu_2-\nu_1/2}{\nu_2} & \textrm{if} & \nu_1\leq 2\nu_2\\
0 & \textrm{if} & \nu_1> 2\nu_2,
\end{array}\right.\\[0.3cm]
\lim_{n\to \infty}\Pi_n^{(\nu_1,\nu_2)}(2,1)&=&\lim_{n\to
\infty}\frac{nP(Y_2>\max\{u_{n,1},u_{n,2}\})}{-\log
\varphi(\boldsymbol{\nu})\eta(\boldsymbol{\nu})}=\left\{
\begin{array}{lcl}
\frac{\nu_1/2}{\nu_2} & \textrm{if} & \nu_1\leq 2\nu_2\\
\frac{\nu_2}{\nu_1/2} & \textrm{if} & \nu_1> 2\nu_2.
\end{array}\right.
\end{eqnarray*}}
\end{ex}\pagebreak



\vspace{1cm}\noindent {\Large  \bf {References }}

\refmark Davis, R.A. (1982). Limit laws for the maximum and minimum
of stationary sequences. {\it{Z. Wahrsch. verw. Gebiete}}. {\bf 61},
31-42.

\refmark Ferreira, H. (1994). Multivariate extreme values in
$T$-periodic random sequences under mild oscillation restrictions.
{\it{Stoch. Proc. App.}}. {\bf 49}, 111-125.

\refmark Ferreira, H. (2006). The upcrossings index and the extremal
index. {\it{J. Appl. Prob.}} {\bf{43}}, 927-937.

\refmark Ferreira, H. (2007). Runs of high values and the
upcrossings index for a stationary sequenced. {\it{Proceedings of
the 56th Session of the International Statistical Institute.}}

\refmark Hsing, T., H\"usler, J. and Leadbetter, M.R. (1988). On the
exceedance point process for a stationary sequence {\it{Prob. Th.
Rel. Fields}}. {{\bf 78}}, 97-112.

\refmark Hsing, T. (1989). Extreme value theory for multivariate
stationary sequences. {\it{J. Multivariate Anal.}} {\bf{29}},
274-291.

\refmark H\"usler, J. (1990). Multivariate extreme values in
stationary random sequences. {\it{Stoch. Proc. Appl.}} {\bf 35},
99-108.

\refmark Nandagopalan, S. (1990). Multivariate Extremes and
Estimation of the Extremal Index. PhD Thesis. University of North
Carolina. Chapel Hill.

\refmark Pereira, L. (2002). Valores Extremos Multidimensionais de
Variaveis Dependentes. PhD Thesis. University of Beira Interior.

\refmark Sebasti\~ao, J., Martins, A.P., Pereira, L. and Ferreira. H., (2010). Clustering of upcrossings of high values. {\it{J. Statist. Plann. Inf.}}, 140, 1003-1012

\end{document}